\documentclass[12pt,twoside,a4paper]{amsart}
\usepackage[utf8]{inputenc}
\usepackage[english]{babel}
\usepackage{float} 
\usepackage[margin=.8in]{geometry}
\usepackage{tkz-euclide}
\usepackage{url}
\usepackage{amsmath,amsfonts,amsthm,amssymb,enumerate,bbm}
\usepackage{color}
\usepackage{hyperref}
\usepackage{graphicx}
\usepackage{enumitem}

\usepackage{setspace}
\hypersetup{linkcolor=blue, colorlinks=true, citecolor = blue}

\allowdisplaybreaks

\numberwithin{equation}{section}

\newtheorem{Theorem}{Theorem}[section]
\newtheorem{Lemma}[Theorem]{Lemma}

\newtheorem{Remark}[Theorem]{Remark}

\let\div\relax
\DeclareMathOperator{\div}{\mathrm{div}}

\newcommand{\Id}{\mathbb{I}}
\newcommand{\dd}{\ \mathrm{d}}

\newcommand{\del}{\partial}
\newcommand{\R}{\mathbb{R}}
\newcommand{\mf}{\mathcal{F}}
\newcommand{\ms}{\mathcal{B}}
\newcommand{\vc}{\mathbf}
\newcommand{\vu}{\vc u}
\newcommand{\bS}{\mathbb S}
\newcommand{\bD}{\mathbb D}

\title[Coll. Non-Newt. and Heat cond. Fluids]{A collision result for both non-Newtonian and heat conducting Newtonian compressible fluids} 

\author{\v S\'arka Ne\v casov\'a}
\address{Institute of Mathematics, Czech Academy of Sciences,
\v Zitn\'a 25, 115 67 Praha 1, Czech Republic.}
\email{matus@math.cas.cz}

\author{Florian Oschmann}
\address{Institute of Mathematics, Czech Academy of Sciences,
\v Zitn\'a 25, 115 67 Praha 1, Czech Republic.}
 \email{oschmann@math.cas.cz}
 
\date{\today}

\begin{document}

\begin{abstract}
We generalize the known collision results for a solid in a 3D compressible Newtonian fluid to compressible non-Newtonian ones, and to Newtonian fluids with temperature depending viscosities.
\end{abstract}

\maketitle

\bigskip




\section{Introduction}

We consider the compressible Navier-Stokes equations governing the motion of a fluid in some bounded domain $\Omega \subset \R^3$, where we additionally insert a simply connected compact obstacle $\ms \subset \R^3$. Denoting $\mf = \Omega \setminus \ms$ the fluid's domain, the equations take the form
\begin{align}\label{NSE}
\begin{cases}
\del_t \rho + \div(\rho \vu) = 0 & \text{in } \mf,\\
\del_t(\rho \vu) + \div(\rho \vu \otimes \vu) - \div \bS + \nabla p = \rho \vc f & \text{in } \mf,\\
m \ddot{\vc G}(t) = -\int_{\del \ms} (\bS - p\Id) \vc n \dd \sigma + \int_\ms \rho_\ms \vc f \dd x & \text{in } \mf,\\
\frac{\rm d}{{\rm d} t}(\mathbb J \omega) = -\int_{\del \ms} (x-\vc G) \times (\bS - p\Id)\vc n \dd \sigma + \int_\ms (x-\vc G) \times \rho_\ms \vc f \dd x & \text{in } \mf,\\
\vu = \dot{\vc G}(t) + \omega(t) \times (x - \vc G(t)) & \text{on } \del \ms,\\
\vu = 0 & \text{on } \del \Omega,\\
\rho(0,\cdot) = \rho_0, \ (\rho \vu)(0,\cdot) = \vc m_0, \ \vc G(0) = \vc G_0, \ \dot{\vc G}(0) = \vc V_0, \ \omega(0) = \omega_0 & \text{in } \mf(0).
\end{cases}
\end{align}
Here, $\rho$ and $\vu$ denote the fluid's density and velocity, respectively, $p$ is the fluid's pressure given by $p(\rho)=\rho^\gamma$ for some $\gamma>\frac32$, $\bS$ the (viscous) stress tensor, and $\vc f \in L^\infty((0,T) \times \Omega)$ a given external force density. Furthermore, $\rho_\ms$ is the solid's density, $\vc G$ the center of mass of the body $\ms$, $\omega$ its rotational velocity, $m>0$ the object's mass given by
\begin{align*}
m=\int_\ms \rho_\ms \dd x,
\end{align*}
and $\mathbb{J}$ is the inertial tensor (moment of inertia) given by
\begin{align*}
\mathbb{J}=\int_\ms \rho_\ms \big( |x-\vc G|^2 \Id - (x-\vc G)\otimes (x-\vc G) \big) \dd x.
\end{align*}

The question of whether or not a solid body collides with its container has been addressed by several authors. Without claiming completeness, we refer to \cite{GerardVaretHillairet2010, GerardVaretHillairet2012, GerardVaretHillairet2014, GVHW2015, Hillairet2007, JNOR2022} for recent results in this direction. The aim of this short note is to generalize these results to non-Newtonian fluids, as well as to heat conducting fluids with temperature growing viscosities. Lastly, let us also mention the related, though different, work \cite{FeireislNecasova2008}, where the authors considered a so-called $k$- or multi-polar compressible fluid, and showed that collisions do not occur for $k\geq 3$ since the velocity and, accordingly, the density enjoy higher regularity. Together with the no-collision results given in the references above, this can be roughly summarized as ``high regularity forbids collision''.\\ 

\paragraph{\bf Notations:} Lebesgue and Sobolev spaces will be denoted in the usual way as $L^p(\Omega)$ and $W^{1,p}(\Omega)$, respectively. We will also denote them for vector- and matrix-valued functions as in the scalar case, that is, $L^p(\Omega)$ instead of $L^p(\Omega; \R^3)$. The Sobolev space of trace-zero functions will be denoted by $W_0^{1,p}(\Omega)$. For each $\mathbb{A},\mathbb{B} \in \R^{3 \times 3}$, we set the Frobenius inner product $\mathbb{A}: \mathbb{B} = \sum_{i,j=1}^3 A_{ij} B_{ij}$. Further, we define the Frobenius norm by $|\mathbb{A}|^2 = \mathbb{A} : \mathbb{A}$. To lean the notation, we will write $a \lesssim b$ if there is a generic constant $C>0$ which is independent of $a$, $b$, and the variables of interest such that $a \leq C b$. The constant might change its value wherever it occurs. The domains occupied by the solid and fluid at time $t \geq 0$ are denoted by $\ms(t)$ and $\mf(t)=\Omega \setminus \ms(t)$, respectively.



\section{General assumptions}
 
Let us start by making precise the assumptions on the fluid and solid. First, the stress tensor $\bS$ will depend on the symmetrized velocity gradient $\bD(\vu)=\frac12 ( \nabla \vu + \nabla^T \vu)$ in a way described in \ref{S1}--\ref{S3} below. Second, we assume that the solid is homogeneous with constant mass density $\rho_\ms>0$. The mass and center of mass of the rigid body are given by
\begin{align*}
    m=\rho_\ms|\ms(0)|, &&  \vc{G}(t)=\frac{1}{m}\int_{\ms(t)}\rho_\ms x \dd x.
\end{align*}
We will also assume that the solid's mass is independent of time, that is, $m=\rho_\ms |\ms(t)|$ for any $t\geq 0$, leading to the density-independent expression $\vc G(t) = |\ms(t)|^{-1} \int_{\ms(t)} x \dd x$.

\subsection{The stress tensor and uniform bounds}

The crucial part in analyzing collisions is to investigate the form of the stress tensor $\bS$. We will make the following assumptions:
\begin{enumerate}[label=(S\arabic*)]
\item Continuity: $\bS$ is a continuous mapping from $\R_{\rm sym}^{3\times 3}$ to itself depending only on the symmetric gradient $\bD(\vu) = \frac12 (\nabla\vu + \nabla^T \vu) \in \R_{\rm sym}^{3 \times 3}$. \label{S1}
\item Monotonicity: For any $\mathbb M, \mathbb N \in \R_{\rm sym}^{3 \times 3}$, we have $[\bS(\mathbb{M}) - \bS(\mathbb{N})] : (\mathbb{M} - \mathbb{N}) \geq 0$.\label{S2}
\item Growth: There are absolute constants $\delta \geq 0$ and $0<c_0\leq c_1<\infty$ such that for some $p>1$ and all $\mathbb M \in \R_{\rm sym}^{3\times 3}$, we have $c_0 |\mathbb{M}|^p - \delta \leq \bS(\mathbb{M}):\mathbb{M} \leq c_1 |\mathbb{M}|^p$.\label{S3}
\end{enumerate}
We note that classical power-law fluids like $\bS = |\bD(\vu)|^{p-2} \bD(\vu)$, but also so-called activated Euler fluids with $\bS=\max\{|\bD(\vu)|-\delta_0, 0\}|\bD(\vu)|^{-1} \bD(\vu)$ for some $\delta_0>0$ fit into this setting. In contrast to the fact that we do not consider temperature in here, we will give another example of temperature growing viscosities in Section~\ref{ch:4}. Note moreover that condition \ref{S3} implies by duality $\bS \in L^{p'}((0,T)\times \Omega)$ since
\begin{align}\label{Sp}
\begin{split}
\|\bS\|_{L^{p'}((0,T)\times \Omega)} &= \sup_{\|\mathbb{M}\|_{L^p((0,T) \times \Omega)} \leq 1} \int_0^T \int_\Omega \bS : \mathbb{M} \dd x \dd t \\
&\leq c_1 \sup_{\|\mathbb{M}\|_{L^p((0,T) \times \Omega)} \leq 1} \int_0^T \int_\Omega |\mathbb{M}|^p \dd x \dd t \leq c_1.
\end{split}
\end{align}

\begin{Remark}
We remark that the question of \emph{existence} of a weak solution to problem \eqref{NSE} is just solved in some special cases, see \cite{Feireisl2003} for Newtonian fluids, and \cite{FeireislLiaoMalek2015} for a special non-Newtonian fluid with bounded divergence of the velocity. On the other hand, for non-Newtonian \emph{incompressible} fluids, existence is shown in \cite{FeireislHillairetNecasova2008}, and in \cite{Necasova2009} for incompressible heat conducting fluids. In those two existence results, collisions cannot occur due to a high regularity of the velocity, in particular, $p \geq 4$ there.
\end{Remark}

To start analyzing the collision behavior, one first needs uniform bounds on the velocity and density. With a slight abuse of notation, we extend the velocity and density as
\begin{align*}
\rho = \begin{cases}
\rho & \text{in } \mf,\\
\rho_\ms &  \text{in } \ms,
\end{cases} &&
\vu = \begin{cases}
\vu & \text{in } \mf,\\
\dot{\vc G}(t) + \omega(t)\times (x-\vc G(t)) & \text{in } \ms.
\end{cases}
\end{align*}
Noticing that the energy inequality obtained in \cite{FeireislLiaoMalek2015} implies in our case
\begin{align*}
\bigg[ \int_\Omega \frac12 \rho |\vu|^2 + \frac{\rho^\gamma}{\gamma-1} \dd x \bigg]_{t=0}^{t=\tau} + \int_0^\tau \int_\Omega \bS : \bD(\vu) \dd x \dd t \leq \int_0^\tau \int_\Omega \rho \vc f \cdot \vu \dd x \dd t
\end{align*}
for almost any $\tau \in [0,T]$, an immediate consequence is the uniform estimate
\begin{align}\label{UnifBds}
\gamma>\frac32, \quad \|\rho\|_{L^\infty(0,T;L^\gamma(\mf(\cdot)))}^\gamma + \|\vu\|_{L^p(0,T;W_0^{1,p}(\Omega))}^p + \|\rho |\vu|^2\|_{L^\infty(0,T;L^1(\Omega))} \lesssim E_0 + 1.
\end{align}
Here, $E_0$ is the initial energy of the system given by
\begin{align}\label{eq:E0}
E_0 = \int_{\mf(0)} \frac{|\vc m_0|^2}{2\rho_0} + \frac{\rho_0^\gamma}{\gamma-1} \dd x + \frac{m}{2} |\vc V_0|^2 + \frac12 \mathbb{J}(0)\omega_0 \cdot \omega_0,
\end{align}
and the ``$+1$'' on the right hand-side of \eqref{UnifBds} is a sole remainder of the force $\vc f$ on the right hand-side of \eqref{NSE}$_2$. Moreover, the implicit constant appearing in \eqref{UnifBds} is independent of the mass $m$ and the final time $T$.

We remark that such bounds also hold true for other models of non-Newtonian fluids such as dissipative (measure-valued) solutions, see \cite{AbbatielloFeireisl2020}. However, the additional Reynolds stress appearing in the momentum equation for such type of solutions is not regular enough for our purposes, in particular, we need to work with weak solutions rather than dissipative ones. Since the present work does not focus on existence of weak solutions, for the definition of such we refer the reader to \cite{FeireislLiaoMalek2015}.

\subsection{The solid's shape and main result}

Throughout the paper, we consider a $C^{1,\alpha}$ solid moving vertically over a flat horizontal surface under the influence of gravity. More precisely, we make the following assumptions (see Figure~\ref{fig1} for the main notations):
\begin{enumerate}[label=(A\arabic*)]
\item The source term is provided by the gravitational force ${\vc f}=-g{\vc e}_3$ and $g>0$. \label{a1}
\item The solid moves along and is symmetric to the $x_3$-axis $\{x_1=x_2=0\}$. \label{a2}
\item The only possible collision point is at $x=0 \in \del \Omega$, and the solid's motion is a vertical translation.\label{a3}
\item Near $r=0$, $\partial \Omega$ is flat and horizontal, where $r=\sqrt{x_1^2+x_2^2}$. \label{a4}
\item Near $r=0$, the lower part of $\partial \ms(t)$ is given by \label{a5}
\begin{align*}
    x_3={ h(t)+r^{1+\alpha}},\ r\leq 2r_0\mbox{ for some small enough } r_0>0.
\end{align*}
\item The collision just happens near the flat boundary of $\Omega$: \label{a6}
\[
\inf_{t>0}\mbox{\rm dist} \left(\ms(t),\partial\Omega \setminus [-2r_0, 2r_0]^2\times\{0\}\right)\geq d_0>0.
\]

\end{enumerate} 

By \ref{a2} and \ref{a3}, we may additionally assume that the position of the solid is characterized by its height $h(t)$, in the sense that 
\begin{align*}
\vc G(t) = \vc G(0)+(h(t)-h(0)) \vc e_3 \quad \text{and} \quad \ms(t)=\ms(0)+(h(t)-h(0)) \vc e_3.
\end{align*}
Note especially that this means that the solid rotates at most around the $x_3$-axis, and so $\omega(t) = \pm |\omega(t)| \vc e_3$. This assumption can be made rigorous for Newtonian incompressible fluids and symmetric initial data in 2D, see \cite{GerardVaretHillairet2010}.\\

\begin{figure}
\centering
\begin{tikzpicture}[scale=.7]
\draw[->] (-6,0) -- (6,0);
\draw[->] (0,0) -- (0,7);
\node at (6,0) [anchor=north] {$r$};
\node at (0,7) [anchor=east] {$x_3$};
\draw (-5,0) rectangle (5,6);
\node at (-6.1,4) [anchor=west] {$\del \Omega$};
\draw[->] (6,4.5)--(6,3.5);
\node at (6,4) [anchor=west] {$-g\vc e_3$};
\draw[black, thick] plot [smooth cycle] coordinates {(0,1) (1,1.5) (2,3) (3,5) (0,5.5) (-3,5) (-2,3) (-1,1.5)};
\draw[<->] (0,0) -- (0,1);
\node at (-1,3) {$\ms$};
\node at (-3,1) {$\mf$};
\node at (-.15,.5) [anchor=west] {$h$};
\node[fill=white] at (2,2) [anchor=west] {$x_3=h+r^{1+\alpha}$};
\draw[->] (2,2) -- (.75,2) -- (.75,1.4);
\draw[dashed] (-1,0) -- (-1,1.5);
\draw[dashed] (1,0) -- (1,1.5);
\node at (1,0) [anchor=north] {$2r_0$};
\draw[dotted] (-.5,0) -- (-.5,1.1);
\draw[dotted] (.5,0) -- (.5,1.1);
\node at (-.5,0) [anchor=north] {$-r_0$};
\end{tikzpicture}
\caption{The body $\ms$ and fluid $\mf$ in the container $\Omega$}
\label{fig1}
\end{figure}

Our main result regarding collision now reads as follows:

\begin{Theorem} \label{theo1}

Let $\gamma>\frac{3}{2}$, $2 \leq p < 3$, $0<\alpha\leq 1$, and $\Omega,\ \ms\subset\mathbb R^3$ be bounded domains of class $C^{1,\alpha}$.
Let $(\rho, \vu, {\vc G})$ be a weak solution to \eqref{NSE} enjoying the bounds \eqref{UnifBds}, let $\bS$ comply with \ref{S1}--\ref{S3}, and assume that \ref{a1}--\ref{a6} are fulfilled. If the solid's mass is large enough, and its initial vertical and rotational velocities are small enough, then the solid touches $\partial \Omega$ in finite time provided
\begin{align}\label{jedna}
\begin{split}
&\alpha<\min \bigg\{ \frac{3-p}{2p-1}, \frac{3(4p\gamma - 3p - 6\gamma)}{p\gamma + 3p + 6\gamma} \bigg\} \quad \text{with}\\
    &\frac32< \gamma \leq 3, \ \frac{6\gamma}{4\gamma-3} < p < 3, \quad \text{or} \quad \gamma>3, \ 2 \leq p < 3. 
    \end{split}
\end{align}
\end{Theorem}

\begin{Remark}
The terms ``large enough'' and ``small enough'' should be interpreted in such a way that inequality \eqref{finalIneq} is satisfied. More precisely, for some constant $C_0>0$ which is independent of $m$ and $T$, we ensure collision provided
\begin{align*}
    C_0 \max\{m^{-1/2}, m^{-3/2} \} \bigg(1 + E_0^{\frac12 + \frac{1}{\gamma} + \frac1p} \bigg) < 1.
\end{align*}
\end{Remark}

\begin{Remark}\label{rem1}
Let us mention a few facts about the above constraints. First, the two expressions inside the minimum stem, as one shall expect, from estimating the diffusive and convective part, respectively.

Second, the restriction $p<3$ is due to the diffusive part, see the estimate of $I_4$ in Section~\ref{sec:PfThm}. Moreover, the requirement $p\geq 2$ stems from the convective term, since we need to estimate the square of the velocity in time. Thus, our result as stated above is just valid for shear-thickening fluids. Omitting convection, Theorem~\ref{theo1} still holds provided
\begin{align}\label{aa1}
\gamma>\frac32, \ \frac{\gamma}{\gamma-1}<p<3, \ \alpha<\min \bigg\{ \frac{3-p}{2p-1}, \frac{9(p\gamma - p - \gamma)}{2p\gamma + 3p + 3\gamma} \bigg\},
\end{align}
hence also allowing for shear-shinning fluids if $\gamma>2$.

Third, the first condition on $p$ and $\gamma$ in \eqref{jedna} can be equivalently stated as $\frac{3p}{4p-6} < \gamma \leq 3$, $2<p<3$.

Lastly, the first fraction inside the minimum in \eqref{jedna} wins precisely if $\gamma \geq \frac{3p}{5p-9}$, and in \eqref{aa1} if $\gamma \geq \frac{3p}{4p-6}$. This seems to be optimal in the sense that for $p=2$, $\alpha=\frac13$ is a ``borderline value'' for the incompressible case, which would (loosely speaking) correspond to $\gamma=\infty$ (see \cite[Section~3.1]{GerardVaretHillairet2012} for details). Moreover, the assumptions in \eqref{jedna} coincide with the requirements on $\alpha$ and $\gamma$ made in \cite{JNOR2022}, where the compressible Newtonian case (corresponding to $p=2$) was considered.
\end{Remark}

\begin{Remark}
As will be immediate from the calculations, the same conclusion holds for non-Newtonian heat conducting fluids such that the assumptions \ref{S1}--\ref{S3} are replaced by
\begin{enumerate}[label=(S\arabic*')]
\item Continuity: $\bS$ is a continuous mapping from $(0,\infty) \times \R_{\rm sym}^{3\times 3}$ to $\R_{\rm sym}^{3 \times 3}$ depending continuously on the temperature $\vartheta>0$ and the symmetric gradient $\bD(\vu) = \frac12 (\nabla\vu + \nabla^T \vu) \in \R_{\rm sym}^{3 \times 3}$.
\item Monotonicity: For any $\vartheta \in (0,\infty)$ and any $\mathbb M, \mathbb N \in \R_{\rm sym}^{3 \times 3}$, we have $[\bS(\vartheta, \mathbb{M}) - \bS(\vartheta, \mathbb{N})] : (\mathbb{M} - \mathbb{N}) \geq 0$.
\item Growth: There are absolute constants $\delta \geq 0$ and $0<c_0\leq c_1<\infty$ such that for some $p>1$, all $\vartheta>0$, and all $\mathbb M \in \R_{\rm sym}^{3\times 3}$, we have $c_0 |\mathbb{M}|^p - \delta \leq \bS(\vartheta, \mathbb{M}):\mathbb{M} \leq c_1 |\mathbb{M}|^p$.
\end{enumerate}
\end{Remark}



\section{Construction of test function and proof of main result}\label{ch3}
In this section, we will define an appropriate test function for the momentum equation that will ensure collision. Let  $(\rho, \vu, {\vc G})$ be a weak solution of \eqref{NSE} satisfying the assumptions \ref{a1}--\ref{a6} in the time interval $(0,T_*)$ before collision and enjoying the bounds \eqref{UnifBds}. From now on we denote $\ms_h=\ms_h(t)=\ms(0)+(h(t)-h(0)){\vc e}_3$ and $\mf_{h}=\mf_h(t)=\Omega\setminus \ms_h(t)$. As mentioned before, the assumption on $\ms(t)$ especially means that the whole configuration is cylindrically symmetric with respect to the $x_3$-axis.\\

Collision can occur if and only if $\lim_{t\to T_*}h(t)=0$. Note further that $\mbox{\rm dist}(\ms_h(t),\partial \Omega)= \min\{h(t),d_0\}$ by assumptions \ref{a2} and \ref{a6}.



\subsection{Test function}\label{sec:testfct}
We will make use of cylindrical coordinates $(r,\theta,x_3)$ with the standard basis $(\vc e_r, \vc e_\theta, \vc e_3)$.
We take the same function as in \cite{GerardVaretHillairet2012} (see also \cite{GerardVaretHillairet2010, GVHW2015}), which is constructed as a function ${\vc w}_h$ associated with the solid particle $\ms_h$ frozen at distance $h$. This function will be defined for $h\in (0,\sup_{t\in [0,T_*)}h(t))$. We see that when $h\to 0$, a cusp arises in $\mf_h$, which is contained in
\begin{align}\label{omegah}
    \Omega_{h,r_0}=\{ x \in \mf_h: 0\leq r<r_0,\ 0\leq  x_3\leq  h+r^{1+\alpha},\ r^2=x_1^2+x_2^2\}.
\end{align}
For the sequel, we fix $h$ as a (small) positive constant and define $\psi(r):= h+r^{1+\alpha}$. Note that the common boundary $\del \Omega_{h,r_0}\cap \del \ms_h$ is precisely given by the set $\{0\leq r\leq r_0,\ x_3=\psi(r)\}$.\\

Let us derive how an appropriate test function inside $\Omega_{h, r_0}$ might look like. In order to get rid of the pressure term, we seek for a function $\vc w_h$ which is divergence-free. Additionally, it shall be rigid on $\ms_h$, and comply with its motion. Thus, our test function shall satisfy
\begin{align*}
\vc w_h |_{\ms_h} = \vc e_3, \quad \vc w_h |_{\del \Omega} = 0, \quad \div \vc w_h = 0,
\end{align*}
hence we choose $\vc w_h = \nabla \times (\phi_h \vc e_\theta)$ for some function $\phi_h(r,x_3)$ to be determined. In cylindrical coordinates, we write $\vc w_h$ as
\begin{align}\label{def:wh}
\vc w_h=-\del_3\phi_h \vc e_r+\frac1r \del_r(r\phi_h)\vc e_3.
\end{align}
The boundary conditions on $\vc w_h$ translate for $\phi_h$ into
\begin{align*}
\del_3 \phi_h(r, 0) = 0, && \frac1r \del_r (r \phi_h)(r, 0)=0,\\
\del_3 \phi_h(r, \psi(r)) = 0 && \frac1r \del_r (r \phi_h)(r, \psi(r)) = 1.
\end{align*}
Further, considering the energy
\begin{align*}
\mathcal{E}=\int_{\mf_h} |\nabla \vc w_h|^2 \dd x
\end{align*}
and anticipating that most of it stems from the vertical motion, that is, from the derivative in $x_3$-direction, we get
\begin{align*}
\mathcal{E}\sim \int_{\mf_h} |\del_3^2 \phi_h|^2 \dd x.
\end{align*}
The Euler-Lagrange equation for the functional $\mathcal{E}$ thus reads $\del_3^4 \phi_h(r, x_3)=0$. A simple calculation now leads to the general form
\begin{align*}
\phi_h(r, x_3) = -\frac32 \bigg( \frac{\kappa_1}{r} - r \bigg) \bigg(\frac{x_3}{\psi(r)} \bigg)^2 + \bigg( \frac{\kappa_1}{r} - r \bigg) \bigg(\frac{x_3}{\psi(r)}\bigg)^3 + \frac{\kappa_2}{r}, \quad \kappa_1, \kappa_2 \in \R.
\end{align*}
In order to get a smooth bounded function $\phi_h$ for all values of $r$ and $x_3$, we choose $\kappa_1=\kappa_2=0$ to infer
\begin{align*}
\phi_h(r, x_3) = \frac{r}{2} \Phi \bigg( \frac{x_3}{\psi(r)} \bigg), \quad \Phi(t) = t^2(3-2t).
\end{align*}
Hence, inside $\Omega_{h, r_0}$, the so constructed function will take advantage of the precise form of the solid. Extending $\phi_h$ in a proper way, we thus can define a proper test function $\vc w_h$ defined in the whole of $\Omega$.

To achieve this, we use a similar method as in \cite{GerardVaretHillairet2010}: define smooth functions $\chi, \eta$ satisfying
\begin{align*}
&\chi=1 \text{ on } (-r_0,r_0)^2 \times (0,r_0), && \chi=0 \text{ on } \Omega\setminus \big( (-2r_0, 2r_0)^2 \times (0, 2r_0) \big), \\ 
&\eta=1 \text{ on } \mathcal{N}_{d_0/2}, && \eta=0 \text{ on } \Omega \setminus \mathcal{N}_{d_0}, 
\end{align*}
where $d_0>0$ is as in \ref{a6}, and $\mathcal{N}_\delta$ is a $\delta$-neighborhood of $\ms(0)$. With a slight abuse of notation for $\phi_h$, set
\begin{align}\label{phih}
\phi_h(r,x_3)= \frac{r}{2} \begin{cases}
    1 & \text{on } \ms_h,\\
    (1-\chi(r,x_3))\eta(r,x_3-h+h(0)) + \chi(r,x_3)\Phi\left(\frac{x_3}{\psi(r)}\right) & \text{on } \Omega \setminus \ms_h,
\end{cases}
\end{align}
and $\vc w_h = \nabla \times (\phi_h \vc e_\theta)$.
Observe that the function $\vc w_h$ satisfies
\begin{align*}
\vc w_h|_{\partial \ms_h}=\vc e_3,\quad \vc w_h|_{\del \Omega}=0,\quad \div\vc w_h=0.
\end{align*}
Indeed, the divergence-free condition is obvious from the definition of $\vc w_h$. Further, since $\phi_h = r/2$ on $\ms_h$, we have $\vc w_h = \vc e_3$ there. Moreover, by definition of $\chi$ and $\eta$, we have $\phi_h=0$ on $\del \Omega \setminus \big( (-2 r_0, 2 r_0)^2 \times \{ 0 \} \big)$ as long as $r_0$ and $h$ are so small that $h+r_0^{1+\alpha} \leq d_0 < r_0$. Lastly, $\phi_h=0$ on $\del \Omega \cap (-r_0, r_0)^2 \times \{0\}$ by definition of $\chi$ and $\Phi(0)=0$, and in the annulus $\big( (-2r_0, 2r_0)^2 \setminus (-r_0, r_0)^2 \big) \times \{ 0 \}$ we use also $\eta(r,h(0))=0$ for $r>\mathfrak{r}_0$ for some $\mathfrak{r}_0\in (d_0, r_0)$ to finally conclude $\vc w_h|_{\del\Omega}=0$, provided $h$ is sufficiently close to zero.

We summarize further properties in the following Lemma, the proof of which is given in \cite[Lemma~3.1]{JNOR2022}:
\begin{Lemma}\label{BdsWh}
$\vc w_h\in C_c^\infty(\Omega)$ and
\begin{align}\label{est1:wh}
\|\del_h\vc w_h\|_{L^\infty(\Omega\setminus \Omega_{h,r_0})} + \|\vc w_h\|_{W^{1,\infty}(\Omega\setminus \Omega_{h,r_0})}\lesssim 1.
\end{align}
Moreover,
\begin{align*}
\|\vc w_h\|_{L^q(\Omega_{h,r_0})} &\lesssim 1 \text{ for any } q<1+\frac3\alpha,\\
\|\del_h \vc w_h\|_{L^q(\Omega_{h,r_0})} + \|\nabla \vc w_h\|_{L^q(\Omega_{h,r_0})} &\lesssim 1 \text{ for any } q<\frac{3+\alpha}{1+2\alpha}.
\end{align*}
\end{Lemma}

\begin{Remark}
The condition $\alpha (q-1)<3$ coming from $\vc w_h$ is consistent with the results of \cite{Starovoitov2003}, where the author showed that collision is forbidden as long as $\alpha(q-1) \geq 3$. Especially, for shapes of class $C^{1,1}$ like balls, this states that no collision can occur as long as $q \geq 4$, which fits the assumptions made in \cite{FeireislHillairetNecasova2008} and \cite{Necasova2009}. Moreover, the difference $q - \frac{2+\alpha}{1+2\alpha}$ occurs in the incompressible two-dimensional setting in \cite[Theorem~3.2]{FilippasTersenov2021} as an optimal value for the solid to move vertically; our fraction $\frac{3+\alpha}{1+2\alpha}$ thus seems like a three-dimensional counterpart to this.
\end{Remark}



\subsection{Estimates near the collision -- Proof of Theorem \ref{theo1}}\label{sec:PfThm}
Let $0<T<T_*$ and let $\zeta\in C^1_c([0,T))$ with $0\leq \zeta\leq 1$, $\zeta'\leq 0$, and $\zeta=1$ near $t=0$. We take $\zeta(t){\vc w}_{h(t)}$ as test function in the weak formulation of the momentum equation \eqref{NSE}$_2$ with right-hand side ${\vc f}=-g{\vc e}_3$, $g>0$. Recalling $\div\vc w_h=0$ and $\del_t \vc w_{h(t)}=\dot{h}(t) \del_h\vc w_{h(t)}$, we get
\begin{align} \label{e1}
\begin{split}
&\int_0^T\zeta\int_\Omega \rho \vu\otimes \vu:\bD({\vc w}_h) \dd x \dd t + \int^T_0\zeta'\int_\Omega \rho \vu\cdot {\vc w}_h \dd x \dd t\\
&\quad + \int^T_0\zeta \dot{h}\int_\Omega \rho \vu\cdot {\del_h \vc w}_h \dd x \dd t - \int^T_0\zeta\int_\Omega {\bS} :\bD({ \vc w}_h) \dd x \dd t\\
= &\int_0^T\zeta \int_\Omega \rho g \vc e_3\cdot\vc w_h \dd x \dd t - \int_\Omega \vc m_0\cdot \vc w_h \dd x\\
= &\int_0^T\zeta \int_{\ms_h} \rho g \vc e_3\cdot\vc w_h \dd x \dd t + \int_0^T\zeta \int_{\mf_h} \rho g \vc e_3\cdot\vc w_h \dd x \dd t - \int_\Omega \vc m_0\cdot \vc w_h \dd x.
\end{split}
\end{align}

Observe that we have $\vc w_h=\vc e_3$ on $\ms_h$, so for a sequence $\zeta_k\to 1$ in $L^1([0,T))$,
\begin{align*}
    &\int^T_0\zeta_k \int_{\ms_h} \rho g{\vc e}_3\cdot {\vc w}_h \dd x \dd t = \int_0^T\zeta_k \int_{\ms_h} \rho_\ms g \dd x \dd t = mg \|\zeta_k\|_{L^1(0,T)} \to mgT.
\end{align*}

In particular, for a proper choice of $\zeta$, it follows that
\begin{multline}\label{mom}
\frac12 mgT \leq \int_0^T\zeta\int_\Omega \rho \vu\otimes \vu:\bD({\vc w}_h) \dd x \dd t + \int^T_0\zeta'\int_\Omega \rho \vu\cdot {\vc w}_h \dd x \dd t + \int^T_0\zeta \dot{h}\int_\Omega \rho \vu\cdot {\del_h \vc w}_h \dd x \dd t\\
\quad - \int^T_0\zeta\int_\Omega {\bS} :\bD({ \vc w}_h) \dd x \dd t
- \int_0^T\zeta \int_{\mf_h} \rho g \vc e_3\cdot\vc w_h \dd x \dd t + \int_\Omega \vc m_0\cdot \vc w_h \dd x
 = \sum_{j=1}^6 I_j.
\end{multline}

We will estimate each $I_j$ separately, and set our focus on the explicit dependence on $T$ and $m$. For the latter purpose, we split each density dependent integral into its fluid and solid part $I_j^f$ and $I_j^\ms$, respectively. The proof follows the same lines as \cite{JNOR2022}, so we will just state the estimates and highlight the differences due to the non-Newtonian setting.\\

$\bullet$ For $I_2^f$, we have
\begin{align*}
    |I_2^f| \lesssim (E_0+1)^{\frac{1}{2\gamma}+\frac12} \quad \text{as long as} \quad \alpha < \frac{3\gamma-3}{\gamma+1}.
\end{align*}

$\bullet$ For $I_2^\ms$, note that due to $\omega(t) = \pm |\omega(t)| \vc e_3$, $\vu |_{\ms_h} = \dot{\vc G}(t) + \omega \times (x - \vc G(t))$, $\vc G(t) = \vc G(0) + (h(t)-h(0))\vc e_3$, $\rho |_{\ms_h} = \rho_\ms>0$, and $\vc w_h |_{\ms_h} = \vc e_3$, we have
\begin{align*}
\int_{\ms_h} \rho \vu \cdot \vc w_h \dd x = \rho_\ms \int_{\ms_h} \big[ \dot{h} \vc e_3 \pm |\omega|\vc e_3 \times (x - \vc G(0) - (h-h(0))\vc e_3) \big] \cdot \vc e_3 \dd x = m \dot{h}.
\end{align*}
Further, from the bounds \eqref{UnifBds}, we infer
\begin{align*}
\sup_{t\in (0,T)} |\dot h|^2 = \sup_{t \in (0,T)} \frac2m \int_{\ms_h} \rho_\ms |\dot h|^2 \dd x \leq \sup_{t \in (0,T)} \frac2m \int_{\ms_h} \rho_\ms |\vc u|^2 \dd x \lesssim \frac2m (E_0+1).
\end{align*}
Hence, by the choice of $\zeta$ such that $|\zeta'|=-\zeta'$ and $\zeta(0)=1+\zeta(T)=1$, we get
\begin{align*}
    |I_2^\ms| \lesssim -\int_0^T \zeta' m |\dot h| \dd t \lesssim \sqrt{m}(E_0 + 1)^\frac12.
\end{align*}

$\bullet$ For $I_3$, observe that $I_3^\ms=0$ due to $\del_h \vc w_h|_{\ms_h} = \del_h \vc e_3=0$. Next, by Sobolev embedding and the bounds \eqref{UnifBds},
\begin{align*}
    \|\vu\|_{L^p(0,T;L^{p^\ast}(\Omega))} \lesssim \|\vu\|_{L^p(0,T;W_0^{1,p}(\Omega))} \lesssim (E_0+1)^\frac1p,
\end{align*}
where we set $p^\ast = 3p/(3-p)$. 
Thus,
\begin{align*}
    |I_3|&=|I_3^f| \leq \int_0^T \zeta |\dot{h}(t)|\, \|\rho\|_{L^\gamma(\mf(\cdot))} \|\vu\|_{L^{p^\ast}(\Omega)}\|\del_h \vc w_h\|_{L^\frac{p^\ast \gamma}{p^\ast (\gamma-1) - \gamma}(\Omega)} \dd t\\
    &\lesssim \|\rho\|_{L^\infty(0,T;L^\gamma(\mf(\cdot)))} \|\vu\|_{L^p(0,T;L^{p^\ast}(\Omega))} \|\del_h \vc w_h\|_{L^\infty(0,T;L^\frac{p^\ast \gamma}{p^\ast (\gamma-1) - \gamma}(\Omega))} \|\zeta \dot{h}\|_{L^{p'}(0,T)} \\
    &\lesssim (E_0+1)^{\frac{1}{\gamma}+\frac1p} \|\dot{h}\|_{L^\infty(0,T)} \|\zeta\|_{L^{p'}(0,T)}\lesssim \sqrt{\frac1m} (E_0+1)^{\frac1\gamma+ \frac1p + \frac12} T^\frac{1}{p'},
\end{align*}
where we have used the estimates \eqref{UnifBds} and Lemma \ref{BdsWh} under the condition
\begin{align*}
    \frac{p^\ast \gamma}{p^\ast (\gamma-1) - \gamma} < \frac{3+\alpha}{1+2\alpha} \Leftrightarrow \alpha < \frac{2p^\ast \gamma - 3p^\ast - 3\gamma}{p^\ast \gamma + p^\ast + \gamma} =
\frac{9(p\gamma - p - \gamma)}{2p\gamma + 3p + 3\gamma}. 
\end{align*}

$\bullet$ Regarding $I_4$, using that $\bS \in L^{p'}((0,T) \times \Omega)$ is bounded by $c_1>0$ (see \eqref{Sp}), we calculate
\begin{align*}
    |I_4| \lesssim \int_0^T \zeta \|\bS\|_{L^{p'}(\Omega)} \|\nabla \vc w_h\|_{L^p(\Omega)} \dd t \leq \|\zeta\|_{L^p(0,T)} \|\bS\|_{L^{p'}((0,T)\times \Omega)} \|\nabla\vc w_h\|_{L^\infty(0,T;L^p(\Omega))} \lesssim T^\frac{1}{p},
\end{align*}
where we have used Lemma \ref{BdsWh} under the condition
\begin{align*}
    p<\frac{3+\alpha}{1+2\alpha} \Leftrightarrow \alpha< \frac{3-p}{2p-1}.
\end{align*}
From here, we get the restriction $p<3$.\\

$\bullet$ For $I_5=I_5^f$,
\begin{align*}
    |I_5| \leq g (E_0+1)^\frac1\gamma T \quad \text{as long as} \quad \alpha< 3(\gamma-1).
\end{align*}

$\bullet$ Similar to $I_2^f$, we have for $I_6^f$
\begin{align*}
    |I_6^f| \leq \|\vc m_0\|_{L^\frac{2\gamma}{\gamma+1}(\mf(0))}\|\vc w_h\|_{L^\infty(0,T;L^\frac{2\gamma}{\gamma-1}(\Omega))} \lesssim \bigg\|\frac{|\vc m_0|^2}{\rho_0}\bigg\|_{L^1(\mf(0))}^\frac12 \|\rho_0\|_{L^\gamma(\mf(0))}^\frac12 \lesssim (E_0+1)^{\frac12+\frac{1}{2\gamma}}.
\end{align*}

$\bullet$ For $I_6^\ms$, since $\vc m_0\cdot \vc w_h |_{\ms_h} = \rho_\ms \vu(0)|_{\ms_h} \cdot \vc e_3 = \rho_\ms \dot h$,
\begin{align*}
    |I_6^s| \lesssim \sqrt{m}(E_0+1)^\frac12.
\end{align*}

$\bullet$ Let us turn to $I_1$. Due to $\vc w_h|_{\ms_h}=\vc e_3$, we see that $I_1^\ms=0$ since $\bD(\vc w_h)=0$ there. Hence, we calculate
\begin{align*}
    |I_1|&= |I_1^f| \lesssim \int_0^T \zeta \|\rho\|_{L^\gamma(\mf_h))} \|\vu\|_{L^{p^\ast}(\Omega)}^2 \|\nabla\vc w_h\|_{L^\frac{p^\ast \gamma}{p^\ast (\gamma-1) - 2\gamma}(\Omega)} \dd t\\
    &\lesssim \|\rho\|_{L^\infty(0,T;L^\gamma(\mf_h))} \|\nabla\vc w_h\|_{L^\infty(0,T;L^\frac{p^\ast \gamma}{p^\ast (\gamma-1) - 2\gamma}(\Omega))} \int_0^T\zeta \|\nabla\vu \|_{L^p(\Omega)}^2 \dd t\\
    &\lesssim (E_0+1)^\frac1\gamma \|\zeta \|_{L^\frac{p}{p-2}(0,T)} \|\nabla\vu\|_{L^p((0,T)\times\Omega)}^2 \lesssim (E_0+1)^{\frac1\gamma+\frac2p} T^{1-\frac2p},
\end{align*}
by using the estimate \eqref{UnifBds} and Lemma \ref{BdsWh} under the condition
\begin{align*}
    \frac{p^\ast \gamma}{p^\ast (\gamma-1) - 2\gamma}<\frac{3+\alpha}{1+2\alpha} \Leftrightarrow \alpha<\frac{2p^\ast \gamma - 3p^\ast - 6\gamma}{p^\ast \gamma + p^\ast + 2\gamma} =
\frac{3(4p \gamma - 3p - 6\gamma)}{p\gamma + 3p + 6\gamma}. 
\end{align*}
Let us emphasize that this term is the only place where the assumption $p\geq 2$ is needed.\\ 

Collecting all the requirements made above, we infer
\begin{align*}
\gamma>\frac32, \ \ 2 \leq p < 3, \ \ p\gamma>p+\gamma, \ \ 4p\gamma>3p+6\gamma,
\end{align*}
which translates into
\begin{align*}
\frac32< \gamma \leq 3, \ \frac{6\gamma}{4\gamma-3} < p < 3, \quad \text{or} \quad \gamma>3, \ 2 \leq p < 3.
\end{align*}

Note further that for any $\gamma\geq \frac32$ and any $\frac{\gamma}{\gamma-1}<p<3$,
\begin{align*}
    \frac{3(4p \gamma - 3p - 6\gamma)}{p\gamma + 3p + 6\gamma}\leq \frac{9(p\gamma - p - \gamma)}{2p\gamma + 3p + 3\gamma} \leq \frac{3\gamma-3}{\gamma+1} \leq 3(\gamma-1) ,
\end{align*}
and that all estimates are independent of the choice of $\zeta$. Hence, we can take a sequence $\zeta_k\to 1$ in $L^r([0,T))$ for some suitable $r > 1$ without changing the bounds obtained. In turn, collecting all estimates above, we finally arrive at
\begin{align*}
    \frac12 mgT \leq C_0 (1+\sqrt{m}+\sqrt{m}^{-1}) \bigg(1 + (E_0+1)^{\frac12+\frac{1}{2\gamma}} + (E_0+1)^\frac12 + (E_0+1)^\frac1p\\
    \quad + (E_0+1)^{\frac1\gamma + \frac1p + \frac12} + g(E_0+1)^\frac1\gamma + (E_0+1)^{\frac1\gamma+\frac2p} \bigg) (1 + T^\frac{1}{p'} + T^\frac{1}{p} + T^{1-\frac2p} + T),
\end{align*}
which, after dividing by $\frac12 mg$ and using Young's inequality on several terms, leads to
\begin{align}\label{infinalIneq}
    T \leq C_0 \max\{m^{-1/2}, m^{-3/2} \} \bigg(1 + E_0^{\frac12 + \frac{1}{\gamma} + \frac1p} \bigg) (1 + T),
\end{align}
where $C_0$ only depends on $p, \gamma, g, \alpha$, the bounds on $\vc w_h$ obtained in Lemma \ref{BdsWh}, and the implicit constant appearing in \eqref{UnifBds}, provided
\begin{align*}
	&\alpha<\min \bigg\{ \frac{3-p}{2p-1}, \frac{3(4p\gamma - 3p - 6\gamma)}{p\gamma + 3p + 6\gamma} \bigg\} \quad \text{with}\\
    &\frac32< \gamma \leq 3, \ \frac{6\gamma}{4\gamma-3} < p < 3, \quad \text{or} \quad \gamma>3, \ 2 \leq p < 3.
\end{align*}

Recalling the definition of $E_0$ from \eqref{eq:E0} as
\begin{align*}
    E_0 &= \int_{\mf(0)}\bigg(\frac{1}{2}\frac{|\vc m_0|^2}{\rho_0} +  \frac{\rho_0^\gamma}{\gamma-1} \bigg) \dd x + \frac{m}{2}|\vc V_0|^2 + \frac12 \mathbb J(0) \omega_0\cdot\omega_0,\\
    \mathbb J(0) &= \int_{\ms_0}\rho_\ms\Big(|x-\vc{G}_0|^2{\mathbb I}-(x-\vc{G}_0)\otimes (x-\vc{G}_0)\Big) \dd x,
\end{align*}

we see that collision can occur only if the solid's mass in \eqref{infinalIneq} is large enough, meaning in fact it's density is very high. Since $E_0$ depends on the solid's mass, we require the solid initially to have low vertical and rotational speed. More precisely, choosing $\vc V_0$ and $\omega_0$ such that $|\vc V_0|, |\omega_0|=\mathcal O(m^{-\frac12})$, and choosing $m$ high enough such that
\begin{align}\label{finalIneq}
    C_0 \max\{m^{-1/2}, m^{-3/2} \} \bigg(1 + E_0^{\frac12 + \frac{1}{\gamma} + \frac1p} \bigg) < 1,
\end{align}
the solid touches the boundary of $\Omega$ in finite time, finishing the proof of Theorem \ref{theo1}.

\begin{Remark}
We see that if, by change, the constant $C_0<1$ small enough, then we can get rid of the assumption on the smallness of $\vc V_0$ and $\omega_0$ by also choosing $m<1$. Indeed, in this case $\max\{m^{-1/2}, m^{-3/2} \} = m^{-3/2}$ and $E_0 \lesssim 1$. Hence, for appropriate values $m<1$ and $C_0 m^{-3/2} < 1$, inequality \eqref{finalIneq} can still be valid.
\end{Remark}



\section{Newtonian flow with temperature-growing viscosity}\label{ch:4}
This section is devoted to investigate a different model for viscosity that does not fit into the assumptions \ref{S1}--\ref{S3}. More precisely, let
\begin{align}\label{Snsf}
\bS=2\mu(\vartheta) \bigg( \bD(\vu) - \frac13 \div \vu \Id \bigg) +\eta(\vartheta) \div \vu \Id,
\end{align}
where the viscosity coefficients $\mu,\eta$ are assumed to be continuous functions on $(0,\infty)$, $\mu$ is moreover Lipschitz continuous, and they satisfy
\begin{align*}
1+\vartheta &\lesssim \mu(\vartheta), \quad |\mu'|\lesssim 1, && 0 \leq \eta(\vartheta) \lesssim 1+\vartheta.
\end{align*}
Note that this means we consider a Newtonian fluid with growing viscosities that are \emph{not} uniformly bounded in the temperature variable.

The equations governing the fluid's motion are now given by
\begin{align}\label{NSF}
\begin{cases}
\del_t \rho + \div(\rho \vu) = 0 & \text{in } \mf,\\
\del_t(\rho \vu) + \div(\rho \vu \otimes \vu) - \div \bS + \nabla p(\rho, \vartheta) = \rho \vc f & \text{in } \mf,\\
m \ddot{\vc G}(t) = -\int_{\del \ms} (\bS - p\Id) \vc n \dd \sigma + \int_\ms \rho_\ms \vc f \dd x & \text{in } \mf,\\
\frac{\rm d}{{\rm d} t}(\mathbb J \omega) = -\int_{\del \ms} (x-\vc G) \times (\bS - p\Id)\vc n \dd \sigma + \int_\ms (x-\vc G) \times \rho_\ms \vc f \dd x & \text{in } \mf,\\
\del_t(\rho s) + \div(\rho s \vu) + \div \frac{\vc q}{\vartheta} = \frac{1}{\vartheta} \big( \bS : \nabla \vu - \frac{\vc q \cdot \nabla \vartheta}{\vartheta} \big) & \text{in } \mf,\\
\vu = \dot{\vc G}(t) + \omega(t) \times (x - \vc G(t)) & \text{on } \del \ms,\\
\vu = 0 & \text{on } \del \Omega,\\
\vc q \cdot \vc n = 0 & \text{on } \del \Omega,\\
\rho(0,\cdot)=\rho_0, \ (\rho \vu)(0,\cdot)=\vc m_0, \ \vartheta(0,\cdot)=\vartheta_0, \ \vc G(0)=\vc G_0, \ \dot{\vc G}(0) = \vc V_0, \ \omega(0) = \omega_0 & \text{in } \mf(0),
\end{cases}
\end{align}
where now $p(\rho, \vartheta) = \rho^\gamma + \rho \vartheta + \vartheta^4$, the heat flow vector $\vc q = \vc q(\vartheta, \nabla \vartheta)$ is given by Fourier's law
\begin{align*}
\vc q(\vartheta, \nabla \vartheta) = - \kappa(\vartheta) \nabla \vartheta
\end{align*}
with the heat conductivity coefficient satisfying
\begin{align*}
\kappa(\vartheta) \sim 1+\vartheta^\beta \ \text{for some} \ \beta>1,
\end{align*}
and the specific entropy $s=s(\rho, \vartheta)$ is connected to the pressure $p(\rho, \vartheta)$ and the internal energy $e(\rho, \vartheta)$ of the fluid by Gibbs' relation
\begin{align*}
\vartheta D s = D e + p D \bigg(\frac{1}{\rho}\bigg).
\end{align*}
Note that this relation determines the internal energy and specific entropy as
\begin{align*}
e(\rho, \vartheta) = \frac{\rho^{\gamma-1}}{\gamma-1} + 3 \frac{\vartheta^4}{\rho} + c_v \vartheta, &&
s(\rho, \vartheta) = 4\frac{\vartheta^3}{\rho} + \log \frac{\vartheta^{c_v}}{\rho},
\end{align*}
where $c_v>0$ is the specific heat capacity at constant volume (see, e.g., \cite{FeireislNovotny2009singlim}). Denoting $\vartheta_\ms>0$ the solid's temperature, we extend the temperature similarly to the velocity and density as
\begin{align*}
\vartheta = \begin{cases}
\vartheta & \text{in } \mf,\\
\vartheta_\ms & \text{in } \ms,
\end{cases}
\end{align*}
and we consider the continuity of the heat flux $\vc q(\vartheta, \nabla\vartheta)\cdot \vc n = \vc q(\vartheta_\ms, \nabla \vartheta_\ms) \cdot \vc n$ on $\del \ms$. Moreover, for simplicity we assume that the heat conductivity coefficient of the solid is the same as the fluid's one (this can be generalized, see \cite[Equation~(4.23)]{Brezina2008}).\\

Noticing that the existence proof of Theorem~4.1.6 in \cite{Brezina2008} also works for any $\beta>2$ instead of $\beta=3$, in such case we have the uniform bound
\begin{align*}
\|\vartheta^\frac{\beta}{2}\|_{L^2(0,T;W^{1,2}(\Omega))}^2 \lesssim E_0+1,
\end{align*}
where this time
\begin{align}\label{initNSF}
E_0 = \int_{\mf(0)} \frac{|\vc m_0|^2}{2 \rho_0} + \rho_0 e(\rho_0, \vartheta_0) \dd x + \frac{m}{2} |\vc V_0|^2 + \mathbb{J}(0)\omega_0 \cdot \omega_0.
\end{align}
Thanks to Sobolev embedding, this yields
\begin{align*}
\vartheta^\frac{\beta}{2} \in L^2(0,T;L^6(\Omega)), \quad \text{that is,} \quad \vartheta \in L^{\beta}(0,T;L^{3\beta}(\Omega)),
\end{align*}
in turn,
\begin{align}\label{bdTheta1}
\|\vartheta\|_{L^\beta(0,T;L^{3\beta}(\Omega))}^\beta \lesssim E_0+1.
\end{align}
Accordingly, the estimate for the stress tensor in $I_4$ changes into
\begin{align*}
|I_4| &= \left| \int_0^T \zeta \int_\Omega \bS:\nabla \vc w_h \dd x \dd t \right| \lesssim \int_0^T \zeta \|\vartheta\|_{L^{3\beta}(\Omega)} \|\nabla \vu\|_{L^2(\Omega)} \|\nabla \vc w_h\|_{L^\frac{6\beta}{3\beta-2}(\Omega)} \dd t \\
&\leq \|\zeta\|_{L^\frac{2\beta}{\beta-2}(0,T)} \|\vartheta\|_{L^\beta(0,T;L^{3\beta}(\Omega))} \|\nabla \vu\|_{L^2((0,T)\times \Omega)} \|\nabla\vc w_h\|_{L^\infty(0,T;L^\frac{\beta}{3\beta-2}(\Omega))}\\
    &\lesssim (E_0+1)^{\frac{1}{\beta} + \frac12} T^{\frac12-\frac{1}{\beta}},
\end{align*}
provided
\begin{align*}
\frac{6\beta}{3\beta-2} < \frac{3+\alpha}{1+2\alpha} \Leftrightarrow \alpha < \frac{3(\beta-2)}{9\beta+2},
\end{align*}
while all the other estimates stay the same. Hence, repeating the arguments from Section~\ref{ch3}, we can state the following
\begin{Theorem}
Let $\gamma>3$, $\beta>2$, $0<\alpha\leq 1$, and $\Omega,\ \ms\subset\mathbb R^3$ be bounded domains of class $C^{1,\alpha}$.
Let $(\rho, \vartheta, \vu, {\vc G})$ be a weak solution to \eqref{NSF} enjoying the bounds \eqref{UnifBds} and \eqref{bdTheta1}, with the initial energy given by \eqref{initNSF}. Moreover, let $\bS$ be given by \eqref{Snsf}, and assume that \ref{a1}--\ref{a6} are satisfied. If the solid's mass is large enough, and its initial vertical and rotational velocities are small enough such that inequality \eqref{finalIneq} is fulfilled, then the solid touches $\partial \Omega$ in finite time provided
\begin{align*}
\alpha < \bigg\{ \frac{3(\gamma-3)}{4\gamma+3}, \frac{3(\beta-2)}{9\beta+2} \bigg\}.
\end{align*}
\end{Theorem}

As can be easily seen, the same arguments can be used for temperature-dependent non-Newtonian fluids, provided the stress tensor decomposes like
\begin{align*}
\bS(\vartheta, \mathbb{M}) = \mu(\vartheta) \tilde{\bS}(\mathbb{M}) + \eta(\vartheta) |\div \vu|^{p-2} \div \vu \Id
\end{align*}
for some tensor $\tilde{\bS}$ satisfying \ref{S1}--\ref{S3}, and $\mu, \eta$ are as above. The details are left to the interested reader.

\begin{Remark}
As a matter of fact, all the analyses in this article also hold for the incompressible case, which (roughly speaking) corresponds to $\gamma=\infty$. Thus, collision for this type of heat conducting fluids occurs if $\beta>2$ and $\alpha<\frac{3(\beta-2)}{9\beta+2}$. Also here, for constant temperature corresponding to a perfectly heat conducting fluid, we recover the borderline value $\alpha<\frac13$ in the limit $\beta\to\infty$, see Remark~\ref{rem1}.
\end{Remark}



\section*{Acknowledgement}
{\it \v S. N. and F. O. have been supported by the Czech Science Foundation (GA\v CR) project 22-01591S.  Moreover, \it \v S. N.  has been supported by  Praemium Academi\ae of \v S. Ne\v casov\' a. The Institute of Mathematics, CAS is supported by RVO:67985840.}



\bibliographystyle{siam}

\end{document}